\documentclass[twoside,11pt]{amsart}
\usepackage{amsmath,amssymb,amsthm}

\begin{document}

\setlength{\textwidth}{145mm} \setlength{\textheight}{203mm}
\setlength{\parindent}{0mm} \setlength{\parskip}{2pt plus 2pt}

\frenchspacing


\numberwithin{equation}{section}
\newtheorem{thm}{Theorem}[section]
\newtheorem{lem}[thm]{Lemma}
\newtheorem{prop}[thm]{Proposition}
\newtheorem{cor}[thm]{Corollary}
\newtheorem{probl}[thm]{Problem}

\newtheorem{defn}{Definition}[section]
\newtheorem{rem}{Remark}[section]
\newtheorem{exa}{Example}



\newcommand{\X}{\mathfrak{X}}
\newcommand{\B}{\mathcal{B}}
\newcommand{\s}{\mathfrak{S}}
\newcommand{\g}{\mathfrak{g}}
\newcommand{\W}{\mathcal{W}}
\newcommand{\LL}{\mathcal{L}}
\newcommand{\dd}{\mathrm{d}}

\newcommand{\diag}{\mathrm{diag}}
\newcommand{\End}{\mathrm{End}}
\newcommand{\im}{\mathrm{Im}}
\newcommand{\id}{\mathrm{id}}

\newcommand{\ie}{i.e. }
\newfont{\w}{msbm9 scaled\magstep1}
\def\R{\mbox{\w R}}
\newcommand{\norm}[1]{\left\Vert#1\right\Vert ^2}
\newcommand{\nnorm}[1]{\left\Vert#1\right\Vert ^{*2}}
\newcommand{\nN}{\norm{N}}
\newcommand{\nJ}{\norm{\nabla J}}
\newcommand{\nnJ}{\nnorm{\nabla J}}
\newcommand{\tr}{{\rm tr}}

\newcommand{\thmref}[1]{Theorem~\ref{#1}}
\newcommand{\corref}[1]{Corollary~\ref{#1}}
\newcommand{\propref}[1]{Proposition~\ref{#1}}
\newcommand{\secref}[1]{\S\ref{#1}}
\newcommand{\lemref}[1]{Lemma~\ref{#1}}
\newcommand{\dfnref}[1]{Definition~\ref{#1}}

\frenchspacing


\title[On the $B$-connection on quasi-K\"ahler manifolds with Norden
metric]{On the geometry of the $B$-connection on quasi-K\"ahler
manifolds with Norden metric}

\author{Dimitar Mekerov}

\maketitle

{\small
\textbf{Abstract} \\
The $B$-connection on almost complex manifolds with Norden metric
is an analogue of the first canonical connection of Lihnerovich in
Hermitian geometry. In the present paper it is considered a
$B$-connection in the class of the quasi-K\"ahler manifold with
Norden metric. Some necessary and sufficient conditions are
derived for the corresponding curvature tensor to be K\"ahlerian.
Curvature properties for this connection are obtained. Conditions
are given for the considered manifolds to be isotropic-K\"ahler.}


\section{Preliminaries}

Let $(M,J,g)$ be a $2n$-dimensional \emph{almost complex manifold
with Norden metric}, \ie $M$ is a differentiable manifold with an
almost complex structure $J$ and a metric $g$ such that
\begin{equation}\label{1.1}
J^2x=-x, \qquad g(Jx,Jy)=-g(x,y)
\end{equation}
for arbitrary $x$, $y$ of the algebra $\X(M)$ on the smooth vector
fields on $M$.

The associated metric $\tilde{g}$ of $g$ on $M$ is defined by
$\tilde{g}(x,y)=g(x,Jy)$. Both metrics are necessarily of
signature $(n,n)$. The manifold $(M,J,\tilde{g})$ is an almost
complex manifold with Norden metric, too.

Further, $x$, $y$, $z$, $w$ will stand for arbitrary elements of
$\X(M)$.

A classification of the almost complex manifolds with Norden
metric is given in \cite{GaBo}. This classification is made with
respect to the tensor field $F$ of type (0,3) defined by
\begin{equation}\label{1.2}
F(x,y,z)=g\bigl( \left( \nabla_x J \right)y,z\bigr),
\end{equation}
where $\nabla$ is the Levi-Civita connection of $g$. The tensor
$F$ has the following properties
\begin{equation}\label{1.3}
F(x,y,z)=F(x,z,y)=F(x,Jy,Jz).
\end{equation}

Among the basic classes $\W_1$, $\W_2$, $\W_3$ of this
classification, the almost complex structure is nonintegrable only
in the class $\W_3$. This is the class of the so-called
\emph{quasi-K\"ahler manifolds with Norden metric}, which we call
briefly \emph{$\W_3$-manifolds}. This class is characterized by
the condition
\begin{equation}\label{1.4}
\mathop{\s} \limits_{x,y,z} F(x,y,z)=0,
\end{equation}
where $\s $ is the cyclic sum over three arguments. The special
class $\W_0$ of the \emph{K\"ahler manifolds with Norden metric}
belonging to any other class is determined by the condition
$F(x,y,z)=0$.

Let $R$ be the curvature tensor of $\nabla$, \ie $R(x,y)z=\nabla_x
\bigl(\nabla_y z\bigr) - \nabla_y \bigl(\nabla_x z\bigr) -
\nabla_{[x,y]}z$. The corresponding tensor of type $(0,4)$ is
determined by $R(x,y,z,w)$ $=g(R(x,y)z,w)$.

The following Ricci identity for almost complex manifolds with
Norden metric is known
\begin{equation}\label{1.5}
\bigl(\nabla_x F\bigr)(y,z,w)-\bigl(\nabla_y
F\bigr)(x,z,w)=R(x,y,Jz,w) - R(x,y,z,Jw).
\end{equation}

The components of the inverse matrix of $g$ are denoted by
$g^{ij}$ with respect to the basis $\{e_i\}$ of the tangent space
$T_pM$ of $M$ at a point $p\in M$.

The \emph{square norm of $\nabla J$} is defined by
\begin{equation}\label{1.6}
    \norm{\nabla J}=g^{ij}g^{ks}
    g\bigl(\left(\nabla_{e_i} J\right)e_k,\left(\nabla_{e_j}
    J\right)e_s\bigr).
\end{equation}
In \cite{MekMan-1} the following equation is proved for a
$\W_3$-manifold
\begin{equation}\label{1.7}
    \norm{\nabla J}=-2g^{ij}g^{ks}
    g\bigl(\left(\nabla_{e_i} J\right)e_k,\left(\nabla_{e_s}
    J\right)e_j\bigr).
\end{equation}

An almost complex manifold with Norden metric $(M,J,g)$ is
K\"ahlerian iff $\nabla J=0$. It is clear that we have
$\norm{\nabla J}=0$ for such a manifold, but the inverse one is
not always true. An almost complex manifold with Norden metric
with $\norm{\nabla J}=0$ is called an \emph{isotropic-K\"ahlerian}
in \cite{MekMan-1}.

The Ricci tensor $\rho$ for the curvature tensor $R$ and the
scalar curvature $\tau$ for $R$ are defined respectively by
\begin{equation}\label{1.8}
    \rho(x,y)=g^{ij}R(e_i,x,y,e_j),\qquad
    \tau=g^{ij}\rho(e_i,e_j),
\end{equation}
and their associated quantities $\rho^*$ and $\tau^*$ are
determined respectively by
\begin{equation}\label{1.9}
    \rho^*(x,y)=g^{ij}R(e_i,x,y,Je_j),\qquad
    \tau^*=g^{ij}\rho(e_i,Je_j).
\end{equation}
Similarly, the Ricci tensor and the scalar curvature are
determined for each \emph{curvature-like tensor (curvature
tensor)} $L$, \ie
 for the tensor $L$ with the following properties:
\begin{equation}\label{1.10}
    L(x,y,z,w)=-L(y,x,z,w)=-L(x,y,w,z),
\end{equation}
\begin{equation}\label{1.11}
    \mathop{\s} \limits_{x,y,z} L(x,y,z,w)=0 \qquad \text{(first Bianchi
    identity)}.
\end{equation}
A curvature-like tensor is called a \emph{K\"ahler tensor} if it
has the property
\begin{equation}\label{1.12}
    L(x,y,Jz,Jw)=-L(x,y,z,w).
\end{equation}


\section{The $B$-connection on $\W_3$-manifolds}

A linear connection $D$ on an almost complex manifold with Norden
metric $(M,J,g)$ preserving $J$ and $g$, \ie $DJ=Dg=0$, is called
a \emph{natural connection} \cite{GaMi}. In \cite{GaGrMi}, on
$(M,J,g)\in\W_1$ a natural connection is considered defined by
\begin{equation}\label{2.1}
    D_x y=\nabla_x y +\frac{1}{2}\bigl(\nabla_x J\bigr)Jy,
\end{equation}
which is called a \emph{$B$-connection}. This connection is known
in Hermitian geometry as the \emph{first canonical connection of
Lihnerovich}.

Let $T$ be a torsion tensor of the $B$-connection $D$ determined
on $(M,J,g)$ by \eqref{2.1}. Because of the symmetry of $\nabla$,
from \eqref{2.1} we have $T(x,y)=\frac{1}{2}\{\bigl(\nabla_x
J\bigr)Jy- \bigl(\nabla_y J\bigr)Jx\}$. Then, having in mind
\eqref{2.1}, we obtain
\[
    T(x,y,z)=g(T(x,y),z)=\frac{1}{2}\{F(x,Jy,z)-F(y,Jx,z)\}.
\]
Let us substitute $z$ by $Jz$. Then, because of \eqref{1.3}, the
last equality implies
\begin{equation}\label{2.2}
    T(x,y,Jz)=\frac{1}{2}\{F(x,y,z)-F(y,x,z)\}.
\end{equation}

We consider the $B$-connection $D$ determined by \eqref{2.1} on a
$\W_3$-manifold $(M,J,g)$. Then, because of \eqref{1.4}, from
\eqref{2.2} we have
\begin{equation}\label{2.3}
    \mathop{\s} \limits_{x,y,z} T(x,y,Jz)=0.
\end{equation}
Therefore, the proposition is valid.
\begin{prop}\label{prop-2.1}
    The torsion tensor $T$ of the $B$-connection $D$ on a
$\W_3$-manifold $(M,J,g)$ satisfies the identity \eqref{2.3}.
\end{prop}

Let $Q$ be the tensor field determined by
\begin{equation}\label{2.4}
    Q(y,z)=\frac{1}{2}\bigl(\nabla_y J\bigr)Jz.
\end{equation}
Having in mind \eqref{1.2}, for the corresponding tensor field of
type (0,3) we have
\begin{equation}\label{2.5}
    Q(y,z,w)=\frac{1}{2}F(y,Jz,w).
\end{equation}
Because of the properties \eqref{1.3}, \eqref{2.5} implies
$Q(y,z,w)=-Q(y,w,z)$.

Let $K$ be the curvature tensor of the $B$-connection $D$, \ie
 $K(x,y)z=$\\ $D_x \bigl(D_y z\bigr) - D_y \bigl(D_x z\bigr) -
D_{[x,y]}z$. Then, according to \eqref{2.1} and \eqref{1.5}, for
the corresponding tensor of type $(0,4)$ we have
\begin{equation}\label{2.6}
\begin{array}{l}
    K(x,y,z,w)=R(x,y,z,w)+\left(\nabla_x Q\right)(y,z,w)-\left(\nabla_y
    Q\right)(x,z,w)
  \\[4pt]
\phantom{K(x,y,z,w)=}
       +Q\left(x,Q(y,z),w\right)-Q\left(y,Q(x,z),w\right).
\end{array}
\end{equation}

After a covariant differentiation of \eqref{2.5}, a substitution
in \eqref{2.6}, a use of \eqref{1.2}, \eqref{1.3}, \eqref{2.4} and
some  calculations, from \eqref{2.6} we obtain
\[
\begin{array}{l}
    K(x,y,z,w)=R(x,y,z,w)-\frac{1}{2}\left[\left(\nabla_x F\right)(y,z,Jw)-\left(\nabla_y
    F\right)(x,z,Jw)\right]
  \\[4pt]
\phantom{K(x,y,z,w)=}
       -\frac{1}{4}\left[g\bigl(\left(\nabla_y J\right)z,\left(\nabla_x J\right)w\bigr)
       -g\bigl(\left(\nabla_x J\right)z,\left(\nabla_y J\right)w\bigr)\right].
\end{array}
\]

The last equality, having in mind \eqref{1.5}, implies
\begin{equation}\label{2.7}
    K(x,y,z,w)=\frac{1}{4}\left\{2R(x,y,z,w)-2R(x,y,Jz,Jw)+P(x,y,z,w)\right\},
\end{equation}
where $P$ is the tensor determined by
\begin{equation}\label{2.8}
    P(x,y,z,w)=g\bigl(\left(\nabla_x J\right)z,\left(\nabla_y J\right)w\bigr)
       -g\bigl(\left(\nabla_y J\right)z,\left(\nabla_x
       J\right)w\bigr).
\end{equation}

In this way, the theorem is valid.
\begin{thm}\label{thm2.2}
    The curvature tensor $K$ of the $B$-connection $D$ on a
$\W_3$-manifold $(M,J,g)$ has the form \eqref{2.7}.
\end{thm}

From \eqref{2.7} is follows immediately that the property
\eqref{1.10} is valid for $K$. Because of $DJ=Dg=0$, the property
\eqref{1.12} for $K$ is valid, too. Therefore, the property
\eqref{1.11} for $K$ is a necessary and sufficient condition $K$
to be a K\"ahler tensor. Since $R$ satisfies \eqref{1.11}, then
from \eqref{2.7} we obtain immediately the following
\begin{thm}\label{thm2.3}
    The curvature tensor $K$ of the $B$-connection $D$ on a
$\W_3$-manifold $(M,J,g)$ is K\"ahlerian iff
\begin{equation}\label{2.9}
    2\mathop{\s} \limits_{x,y,z} R(x,y,Jz,Jw)=\mathop{\s} \limits_{x,y,z}
    P(x,y,z,w).
\end{equation}
\end{thm}

Let the following condition be valid for the $\W_3$-manifold
$(M,J,g)$

\begin{equation}\label{2.10}
    \mathop{\s} \limits_{x,y,z} R(x,y,Jz,Jw)=0.
\end{equation}

The condition \eqref{2.10} characterizes the class $\LL_2$ of the
almost complex manifolds with Norden metric according to the
classification in \cite{GrDjMe} with respect to the properties of
$R$.

The equality \eqref{2.8} implies immediately the properties
\eqref{1.10} and \eqref{1.12} for $P$. Then, according to
\eqref{2.9} and \eqref{2.10}, we obtain the following
\begin{thm}\label{thm2.4}
    Let $(M,J,g)$ belongs to the class $\W_3\cap\LL_2$. Then the
    curvature tensor $K$ of the $B$-connection $D$ is K\"ahlerian iff
    the tensor $P$ determined by \eqref{2.8} is K\"ahlerian, too.
\end{thm}

Having in mind \eqref{2.7}, the last theorem implies the following
\begin{cor}\label{cor2.5}
    Let the curvature tensor $K$ of the $B$-connection $D$ be
    K\"ah\-ler\-ian on $(M,J,g)\in\W_3\cap\LL_2$. Then the tensor $H$,
    determined by
    \begin{equation}\label{2.11}
        H(x,y,z,w)=R(x,y,z,w)-R(x,y,Jz,Jw)
    \end{equation}
    is a K\"ahler tensor.
\end{cor}


\section{Curvature properties of the connection $D$ in $\W_3\cap\LL_2$}

Let us consider the manifold $(M,J,g)\in\W_3\cap\LL_2$ with
K\"ahler curvature tensor $K$ of the $B$-connection $D$. Then,
according to \thmref{thm2.4} and \corref{cor2.5}, the tensor $P$
and $H$, determined by \eqref{2.8} and \eqref{2.11}, respectively,
are also K\"ahlerian.

Let $\rho(K)$ and $\rho(P)$ be the Ricci tensors for $K$ and $P$,
respectively. Then we obtain immediately from \eqref{2.7}
    \begin{equation}\label{3.1}
        \rho(y,z)-\rho^*(y,Jz)=2\rho(K)(y,z)-\frac{1}{2}\rho(P)(y,z).
    \end{equation}
We denote $\tau^{**}=g^{ij}g^{ks}R(e_i.e_k,Je_s,Je_j)$ and from
\eqref{3.1} we have
    \begin{equation}\label{3.2}
        \tau-\tau^{**}=2\tau(K)-\frac{1}{2}\tau(P),
    \end{equation}
where $\tau(K)$ and $\tau(P)$ are the scalar curvatures for $K$
and $P$, respectively.

It is known from \cite{MekMan-1}, that $\nJ=-2(\tau+\tau^{**})$.
Then \eqref{3.2} implies
    \begin{equation}\label{3.3}
        \tau=\tau(K)-\frac{1}{4}\left(\tau(P)+\nJ\right).
    \end{equation}

From \eqref{2.8} we obtain
\[
    \rho(P)(y,z)=g^{ij}g\bigl(\left(\nabla_{e_i} J\right)z,\left(\nabla_y
    J\right)e_j\bigr),
\]
from where
\[
    \tau(P)=g^{ij}g^{ks}g\bigl(\left(\nabla_{e_i} J\right)e_s,\left(\nabla_{e_k}
    J\right)e_j\bigr).
\]
Hence, applying \eqref{1.7}, we get
    \begin{equation}\label{3.4}
        \tau(P)=-\frac{1}{2}\nJ.
    \end{equation}
From \eqref{3.3} and \eqref{3.4} it follows
    \begin{equation}\label{3.5}
        \tau=\tau(K)-\frac{1}{8}\nJ.
    \end{equation}

The equality \eqref{3.5} implies the following
\begin{prop}\label{prop3.1}
    Let the curvature tensor $K$ of the $B$-connection $D$ be
    K\"ahlerian on $(M,J,g)\in\W_3\cap\LL_2$. Then $(M,J,g)$ is an isotropic-K\"ahler manifold iff
    $\tau=\tau(K)$.
\end{prop}

Let the considered manifold with K\"ahler curvature tensor $K$ of
the $B$-connection $D$ in $\W_3\cap\LL_2$ be 4-dimensional. Since
$H$ is a K\"ahler tensor, then according to \cite{Teof-2}, we have
\begin{equation}\label{3.6}
      H=\nu(H)(\pi_1-\pi_2)+\nu^*(H)\pi_3,
\end{equation}
where $\nu(H)=\frac{\tau(H)}{8}$, $\nu^*(H)=\frac{\tau^*(H)}{8}$
and
\[
\begin{array}{l}
\pi_1(x,y,z,w)=g(y,z)g(x,w)-g(x,z)g(y,w),\\[4pt]
\pi_2(x,y,z,w)=g(y,Jz)g(x,Jw)-g(x,Jz)g(y,Jw),\\[4pt]
\pi_3(x,y,z,w)=-g(y,z)g(x,Jw)+g(x,z)g(y,Jw),\\[4pt]
\phantom{\pi_2(x,y,z,w)=} -g(y,Jz)g(x,w)+g(x,Jz)g(y,w).
\end{array}
\]

From \eqref{3.6}, \eqref{2.7} and \eqref{2.8} follows the next
\begin{prop}\label{prop3.2}
    Let the curvature tensor $K$ of the $B$-connection $D$ be
    K\"ahlerian on a 4-dimensional $(M,J,g)\in\W_3\cap\LL_2$ and $P$ and
    $H$ are determined by \eqref{2.8} and \eqref{2.11}, respectively.
    Then $H$ has the form
    \[
          H=\frac{4\tau(K)-\tau(P)}{16}(\pi_1-\pi_2)+\frac{4\tau^*(K)-\tau^*(P)}{16}\pi_3.
    \]
\end{prop}

\bigskip


\bigskip

\textit{Dimitar Mekerov\\
University of Plovdiv\\
Faculty of Mathematics and Informatics
\\
Department of Geometry\\
236 Bulgaria Blvd.\\
Plovdiv 4003\\
Bulgaria
\\
e-mail: mircho@uni-plovdiv.bg}

\end{document}